\documentclass[conference,a4paper]{IEEEtran}

\usepackage[colorinlistoftodos]{todonotes}

\IEEEoverridecommandlockouts
\usepackage{cite}
\usepackage[numbers]{natbib}

\makeatletter
\renewcommand{\thebibliography}[1]{%
  \section*{\refname}%
  \footnotesize
  \list{\@biblabel{\@arabic\c@enumiv}}%
       {\settowidth\labelwidth{\@biblabel{#1}}%
        \leftmargin\labelwidth
        \advance\leftmargin\labelsep
        \usecounter{enumiv}}%
  \sloppy
  \clubpenalty4000
  \@clubpenalty \clubpenalty
  \widowpenalty4000%
  \sfcode`\.\@m}
\makeatother

\usepackage{amsmath,amssymb,amsfonts}
\usepackage{graphicx}
\usepackage{textcomp}
\usepackage{xcolor}
\usepackage{multirow}
\usepackage{array}
\usepackage{algorithm}
\usepackage{algpseudocode}
\usepackage{caption}

\usepackage{mathtools}

\makeatletter
\DeclareRobustCommand\onedot{\futurelet\@let@token\@onedot}
\def\@onedot{\ifx\@let@token.\else.\null\fi\xspace}

\makeatother

\definecolor{red}{rgb}{1,0,0}
\definecolor{orange}{rgb}{1,.5,0}
\definecolor{darkgreen}{rgb}{0,.5,0}
\definecolor{grey}{rgb}{0.5,0.5,0.5}

\newif\ifcomments
\newif\ifshortmodel
\newif\iflongmodel
\shortmodeltrue

\newcommand{\PH}[1]{}
\newcommand{\ET}[1]{}
\newcommand{\TD}[1]{}

\ifcomments
\renewcommand{\PH}[1]{\textcolor{red}{Philipp: #1}}
\renewcommand{\ET}[1]{\textcolor{orange}{Eric: #1}}
\renewcommand{\TD}[1]{\textcolor{teal}{TODO: #1}}
\fi

\def\alldems{_{d}}
\def\alllines{_{(i,j)}}

\def\genatbus{_{g|\mathcal{B}(g)=i}}
\def\dematbus{_{d|\mathcal{B}(d)=i}}

\def\linesatbus{_{j \in I}}

\def\gen{^{\text{G}}}
\def\genmin{^{\text{G}}}
\def\genmax{^{\text{G}}}

\def\dem{^{\text{D}}}

\def\sol{n}

\def\xattm{x^{(n)}_{\mathrm{a},\overline Z}}

\def\xattnew{$x^{(n+1)}_{\mathrm{a},\overline Z}$}
\def\xatt{$\xattm$}
\def\xattac{x^{(n)}_{\mathrm{LAC},\overline Z}}

\def\hidden{U}

\def\acset{\mathcal{N}^{\text{LAC}}}

\def\objacm{\zeta^{(n)}_{\mathrm{LAC},\overline Z}}
\def\objdcm{\zeta^{(m(n))}_{\mathrm{DC},\overline Z}}
\def\obja{\zeta^{(n)}_{\mathrm{a},\overline Z}}
\def\objabs{\zeta^{(n)}_{\text{abs},\overline Z}}
\def\objrel{\zeta^{(n)}_{\text{rel},\overline Z}}

\def\BibTeX{{\rm B\kern-.05em{\sc i\kern-.025em b}\kern-.08em
    T\kern-.1667em\lower.7ex\hbox{E}\kern-.125emX}}
\begin{document}

\allowdisplaybreaks

\title{Vulnerability Analysis Evaluating Bilevel Optimal Power Flow Approaches for Multiple Load Cases
\thanks{This work is funded by the Deutsche Forschungsgemeinschaft (DFG, German Research Foundation), project number 360352892, priority programme DFG SPP 1984.}
}

\author{
\IEEEauthorblockN{Eric Tönges}
\IEEEauthorblockA{\textit{Sustainable Electrical Energy Systems} \\
\textit{University of Kassel}\\
Kassel, Germany \\
eric.toenges@uni-kassel.de}
\and
\IEEEauthorblockN{Martin Braun}
\IEEEauthorblockA{\textit{Sustainable Electrical Energy Systems} \\
\textit{University of Kassel, Fraunhofer IEE}\\
Kassel, Germany \\
martin.braun@iee.fraunhofer.de}
\and
\IEEEauthorblockN{Philipp Härtel}
\IEEEauthorblockA{\textit{Sustainable Electrical Energy Systems} \\
\textit{University of Kassel, Fraunhofer IEE}\\
Kassel, Germany \\
philipp.haertel@iee.fraunhofer.de}
\and
}

\IEEEoverridecommandlockouts
\maketitle
\IEEEpubidadjcol

\begin{abstract}

This work presents two methodologies to enhance vulnerability assessment in power systems using bilevel attacker-defender network interdiction models. 
First, we introduce a systematic evaluation procedure for comparing different optimal power flow formulations in the lower-level problem. 
We demonstrate the procedure for a comparison of the widely used DC approximation and a linearized AC optimal power flow model. 
Second, we propose a novel scoring methodology to identify and prioritize critical attack vectors across diverse load and generation scenarios.
Both methodologies go beyond traditional worst-case analysis.
Case studies on a SimBench high-voltage test grid show that the DC approach fails to detect a significant portion of critical vulnerabilities.
The scoring methodology further demonstrates the dependency of vulnerabilities on the considered load case and time step, highlighting the importance of assessing multiple scenarios and going beyond worst-case solutions.
The proposed methodologies enhance power system vulnerability assessment and can support the effective development of robust defense strategies for future power systems.

\end{abstract}


\section*{Nomenclature}
\subsection{Sets, indices and mapping functions}
\addcontentsline{toc}{section}{Nomenclature}
\begin{IEEEdescription}[\IEEEusemathlabelsep\IEEEsetlabelwidth{$V_1,V_2,V_3$}]
\item[$d \in D$] Index and set of all demand units
\item[$g \in G$] Index and set of all generation units
\item[$i, j \in I$] Index and set of all buses
\item[$(i,j) \in K$] Index and set of all branches, each $(i,j)$ has a corresponding $(j,i)$
\item[$\mathcal{B}(\cdot)\in I$] Bus of generator $g$ or demand $d$
\item[$\Omega_{(\cdot)}$] Set of lower-level decision variables of the considered problem
\end{IEEEdescription}
\vspace{-0.3em}
\subsection{Parameters}
\addcontentsline{toc}{section}{Nomenclature}
\begin{IEEEdescription}[\IEEEusemathlabelsep\IEEEsetlabelwidth{$V_1,V_2,V_3$}]
\item[$B_{ij}, G_{ij}$] Susceptance and conductance of branches $(i,j)$ and $(j,i)$
\item[$B_{ij}^\text{S}$] Shunt susceptance of branches $(i,j)$ and $(j,i)$ 
\item[$ P_{d}^{\text{D}}$] Input active-power consumption of demand $d$\vspace{0.2em}
\item[$ \overline P_{g}^{\text{G}}$] Max. active-power injection of generator $g$ \vspace{0.1em}
\item[$\underline Q_g^G, \overline Q_g^G$] Min. and max. reactive power of generator $g$\vspace{0.1em}
\item[$\overline S_{ij}$] Max. apparent-power flow of branch $(i,j)$ and $(j,i)$
\item[$\underline V_i, \overline V_i$] Min. and max. voltage magnitude of bus $i$ 
\item[$\overline Z$] Max. number of attacked branches
\item[$\alpha_{d}^{\text{D}}, \alpha_{g}^{\text{G}}$] Active-reactive power ratio of demand~$d$ or generator $g$
\end{IEEEdescription}
\subsection{Decision variables}
\addcontentsline{toc}{section}{Nomenclature}
\begin{IEEEdescription}[\IEEEusemathlabelsep\IEEEsetlabelwidth{$\delta_d^{\text{p,D}(\cdot)}, \delta_g^{\text{p,G}}$]}]
\item[$p_d\dem, q_d\dem$] (Re-)active power of demand $d$ \vspace{0.2em}
\item[$p_g\gen, q_g\gen$] (Re-)active power of generator $g$ \vspace{0.2em}
\item[$p_{ij}^{\text{K}}, q_{ij}^\text{K}$] (Re-)active power injected into branch $(i,j)$ at bus $i$
\item[$v_i, \theta_i$] Voltage magnitude and angle of bus $i$
\item[$z_{ij} \in \{0,1\}$] Binary in-service variable of branch $(i,j)$ (equals 1 if in service, 0 if attacked)
\end{IEEEdescription}
\subsection{Additional definitions for proposed methodologies}
\addcontentsline{toc}{section}{Nomenclature}
\begin{IEEEdescription}[\IEEEusemathlabelsep\IEEEsetlabelwidth{$a \in \{ \text{DC}, \text{LAC}\}$}]
\item[$a \in \{ \text{DC}, \text{LAC}\}$] Index of the modeling approach
\item[$N^{a}$] Max. number of considered CAVs
\item[$\mathcal{N}^a = \{1, ..., N^{a}\}$] Ordered set of natural numbers of CAVs
\item[$n, m \in \mathcal{N}^a$] CAV indices and set of CAV indices
\item[$\mathcal{L}^{a}$] CAV list, storing all identified CAVs
\item[$t \in \{1, ..., T\}$] Index of time step for $T$ time steps
\item[$u, U$] Relative and absolute undetected CAVs
\item[$x \in \mathcal{X}^{\overline Z}$] Index and set of all possible attack vectors containing $\overline Z$ components
\item[\xatt] $n^{\text{th}}$-best attack vector for given $a, \overline Z$
\item[$z_{ij}^{\star, a, \overline Z}$] Optimal value of $z_{ij}$ for given $a, \overline Z$
\item[$\zeta^{(n)}_{a,\overline Z}$] Optimal objective value for the $n^\text{th}$-best CAV with approach $a$
\item[$\mathcal{C}^{\overline Z}(x)$] Appearance counter over all time steps in which attack combination $x$ appears 
\item[$\mathcal{R}^{\overline Z}(x)$] Sum of ranks $n$ of attack combination $x$
\item[$\mathcal{Y}^{\overline Z}(x)$] Sum of objective values of attack combination $x$
\end{IEEEdescription}

\if False
\subsection{Dual decision variables}
\addcontentsline{toc}{section}{Nomenclature}
\begin{IEEEdescription}[\IEEEusemathlabelsep\IEEEsetlabelwidth{$1234567890123456789012345$}]
\item[$\underline\gamma_d^\text{D}, \overline\gamma_d^\text{D}, \underline\gamma_g^\text{G}, \overline\gamma_g^\text{G}, \epsilon_d^\text{D}, \epsilon_g^\text{G}, \eta_{k,l}^{\rightarrow}, \eta_{k,l}^{\leftarrow},$]
\item[$\underline\eta_k, \overline\eta_k, \kappa_k, \lambda_k^\rightarrow, \lambda_k^\leftarrow, \mu_k^\rightarrow, \mu_k^\leftarrow, \underline\nu_d^\text{D},$] \vspace{0.2em}
\item[$\overline\nu_d^\text{D}, \underline\nu_g^\text{G}, \overline\nu_g^\text{G}, \xi_{km}, \underline\rho_n, \overline\rho_n, \sigma_n^p, \sigma_n^q, $] Dual variables associated
\item[$\tau_k, \underline\varphi_d^\text{D}, \overline\varphi_d^\text{D}, \underline\varphi_g^\text{G}, \overline\varphi_g^\text{G}, \chi_k, \psi_d^\text{D,p}, $] with the problem \vspace{0.2em}
\item[$\psi_d^\text{D,q}, \psi_g^\text{G,p}, \psi_g^\text{G,q}, \underline\omega_d^\text{D}, \overline\omega_d^\text{D}, \underline\omega_g^\text{G}, \overline\omega_g^\text{G}$]
\end{IEEEdescription}
\fi

\section{Introduction} \label{sec:introduction}
The electric power system is a core critical infrastructure and crucial for modern societies. 
Beyond reliability against usual disturbances such as random component failures, resilience against so-called ``high-impact low-probability'' (HILP) events is essential for secure power system planning and operation \cite{Braun.2024}.
Incidents such as the 2015 Ukraine blackout \cite{Liang.2017} and the 2025 Cannes blackout \cite{Asaf.2025} highlight the relevance of deliberate adversarial attacks as sources of HILP events.
Identifying and understanding vulnerabilities is therefore necessary for developing effective resilience-enhancing strategies.
\par
\begin{figure}[b]
\centering
\includegraphics[width=8.83cm]{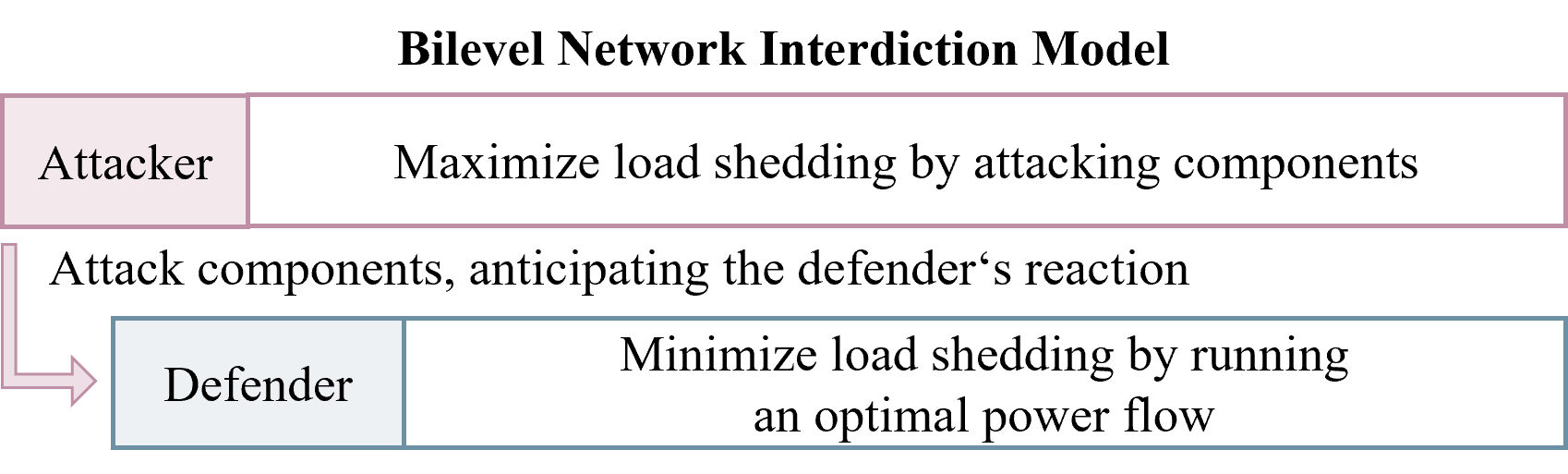}
\caption{General scheme of a BNI model.}
\label{fig:bilevel_scheme}
\end{figure}

Mathematical bilevel optimization, or bilevel network interdiction (BNI) modeling, is a prominent approach for identifying vulnerabilities to adversarial HILP events in power systems \cite{HernandezValencia.2021}. 
In this framework, the upper level represents an attacker seeking to maximize system damage by targeting components while anticipating the grid operator’s response. 
The lower-level model represents the grid operator's efforts to minimize damage based on the upper-level decision, typically through an optimal power flow (OPF) formulation.
Fig.~\ref{fig:bilevel_scheme} depicts the general scheme of a BNI model. 
Since the bilevel model with an exact AC OPF in the lower level is strongly NP-hard due to nonconvexity \cite{Bienstock.2019}, the OPF is often approximated by its linear DC form \cite{Salmeron.2004, Zhao.2013}.
This linearization enables dual reformulations and reduces computational effort.
However, recent studies have proposed several approaches that incorporate an AC OPF in the lower level.
For example, a second-order Taylor approximation is used in \cite{Abedi.2021}, a second-order cone relaxation based on Jabr's inequality \cite{Jabr.2006} is applied in \cite{Wu.2018}, and the authors of \cite{Dandurand.2021, LopezLezama.2017} address the nonconvex AC OPF using branch-and-bound and iterated local search algorithms, respectively. 
Linearized and nonlinear AC OPF formulations were shown to produce more precise worst-case solutions and even different attack vectors compared to DC approaches \cite{Abedi.2021, LopezLezama.2017}, with nonconvex formulations in some cases outperforming AC approaches using second-order cone linearization \cite{Dandurand.2021}.

Most analyses focus on worst-case or predefined attack vectors.
While critical attack vectors (CAVs) beyond the worst case are identified in \cite{Ding.2017}, these are not used to compare different modeling approaches.
In addition, most studies consider only a single grid configuration with fixed load and generation.
Exceptions include \cite{Abedi.2020}, which examines worst-case solutions over multiple time steps, and \cite{Sundar.2024}, which uses a stochastic BNI model.
Two main research gaps remain: First, there is no method to comprehensively assess the quality of different lower-level OPF formulations beyond worst-case scenarios. 
It is thus unclear to what extent linearized models, such as the DC approximation, systematically miss impactful HILP events or whether they simply reorder critical vulnerabilities.
Second, there is no framework combining CAV analysis beyond the worst-case scenario across multiple time steps and load cases.
\par
The main contributions of this work are:
\begin{itemize}
    \item an evaluation procedure for comparing vulnerability assessment approaches with respect to critical attack vectors, moving beyond worst-case analysis and applied to both a linearized AC and a DC bilevel network interdiction model,
    \item a scoring methodology for quantifying the impact of various power flow cases and consistently evaluating critical attack vectors across them, thereby highlighting the benefits of considering multiple cases and attack vectors,
    \item case studies demonstrating that more than 15\,\% of critical attack vectors remain undetected with the DC approach, and that attack vectors with objective scores exceeding 77\,\% of the most critical solution are missed without the proposed scoring methodology.
\end{itemize}
Although demonstrated in the context of BNI modeling, the proposed algorithms and methodologies are broadly applicable to vulnerability assessment approaches and are not limited to bilevel programming.
The results support the selection of appropriate modeling approaches and inform research on resilience-enhancing strategies for power systems.
\par
The remainder is organized as follows. 
Sec.~\ref{sec:formulations} presents the model formulations and describes reformulations.
Sec.~\ref{sec:methodologies} introduces the new methodologies described above, which are applied in case studies in Sec.~\ref{sec:case_studies}.
Sec.~\ref{sec:conclusion} concludes the analysis and provides an outlook on future work.

\section{Model formulations and reformulations for bilevel power system vulnerability assessment} \label{sec:formulations}
This section presents the BNI model formulations with an AC OPF and a DC OPF lower level, as applied in Sec.~\ref{sec:case_studies}, and describes the corresponding linearizations and dual reformulations.

\subsection{Nonconvex AC bilevel model} \label{sec:nonconvex}
The upper level of the nonconvex AC BNI model is given in \eqref{eq:ul_obj}--\eqref{eq:att_symmetrie} and the lower level is presented in \eqref{eq:ll_obj}--\eqref{eq:linelim_from}. The set of lower-level decision variables is defined as $\Omega_{\text{AC}} \coloneqq \big\{p\dem_d, p\gen_g, q\dem_d, q\gen_g, v_i, \theta_i, p_{ij}^\text{K}, q_{ij}^\text{K} \big\} $. 
\par\vspace{-0.5em}
\begin{subequations}\label{eq:bilevel_problem}
\renewcommand{\theequation}{\theparentequation.\arabic{equation}} \vspace{-0.5em}
\begin{flalign}
  & \max_{z_{ij}} \textstyle  \sum\alldems (P_d - p_d\dem) \label{eq:ul_obj} \\
  & \text{s.t.} \nonumber \\
  & \textstyle \sum\alllines 0.5(1-z_{ij}) \leq \overline Z, \quad \text{where } z_{ij} \in \{0, 1\} \quad \forall (i,j) \in K \label{eq:att_budget}\\
  & z_{ij} = z_{ji} \quad \forall (i,j) \in K \label{eq:att_symmetrie} \\
  & \min_{\Omega_{\text{AC}}} \textstyle  \sum\alldems (P_d - p_d\dem) \label{eq:ll_obj}\\
  & \text{s.t.} \nonumber \\
  & 0 \leq p_d\dem \leq P_d\dem \quad \forall d \in D \label{eq:demlim}\\
  & 0 \leq p_g\gen \leq \overline P_g\genmax \quad \forall g \in G \label{eq:plim_gen}\\
  & \underline Q_g\genmin \leq q_g\gen \leq \overline Q_g\genmax \quad \forall g \in G \label{eq:qlim_gen}\\
  & q_d\dem = \alpha_d\dem p_d\dem \quad \forall d \in D\label{eq:cosphi_dem}\\
  & q_g\gen \leq \alpha_g\gen p_g\gen \quad \forall g \in G \label{eq:cosphi_gen}\\
  & \underline V_i \leq v_i \leq \overline V_i \quad \forall i \in I \label{eq:volt_limits}\\
  & \textstyle \sum\genatbus p_g\gen - \sum\dematbus p_d\dem - \textstyle \sum\linesatbus p_{ij}^\text{K} = 0 \quad \forall i \in I \label{eq:p_balance}\\
  & \textstyle \sum\genatbus q_g\gen - \sum\dematbus q_d\dem - \textstyle \sum\linesatbus q_{ij}^\text{K} = 0 \quad \forall i \in I \label{eq:q_balance}\\
  & p_{ij}^\text{K} = z_{ij} \Big[G_{ij} v_i^2 - v_i v_j \big(G_{ij} \cos(\theta_i - \theta_j) \nonumber \\ 
  & + B_{ij} \sin(\theta_i - \theta_j)\big)\Big] \quad \forall (i,j) \in K \label{eq:line_p}\\
  & q_{ij}^\text{K} = z_{ij} \Big[-(B_{ij} + 0.5B_{ij}^\text{S})v_i^2 + v_i v_j \big(B_{ij} \cos(\theta_i - \theta_j) \nonumber \\
  & - G_{ij}\sin(\theta_i - \theta_j)\big)\Big] \quad \forall (i,j) \in K \label{eq:line_q}\\
  & (p_{ij}^\text{K})^2 + (q_{ij}^\text{K})^2 \leq (\overline S_{ij})^2 \quad \forall (i,j) \in K \label{eq:linelim_from}
\end{flalign}
\end{subequations} \vspace{-1.2em}
\par
The upper-level objective \eqref{eq:ul_obj} is to maximize active-power load shedding, which is equivalently minimized in the lower-level objective function \eqref{eq:ll_obj}, subject to solving an OPF. 
The attacker's resources are limited by \eqref{eq:att_budget}, where the binary variable $z_{ij}$ indicates whether a branch is in service. The attacker can select a set of up to $\overline Z$ branches to be set out of service. For each attacked branch, $z_{ij}$ and the corresponding $z_{ji}$ both equal $0$, while $z_{ij}=z_{ji}=1$ holds for in-service branches. It is assumed that an attack is always successful. Since each attack on branch $(i,j)$ involves two variables $z_{ij}$ and $z_{ji}$, \eqref{eq:att_budget} includes a factor of $0.5$.
In the lower level, active-power demand is bounded by \eqref{eq:demlim}, and active- and reactive power generation limits are given in \eqref{eq:plim_gen}--\eqref{eq:qlim_gen}.
The ratio between active and reactive power for demand and generation is constrained by \eqref{eq:cosphi_dem}--\eqref{eq:cosphi_gen}, and \eqref{eq:volt_limits} defines the bounds on voltage magnitude.
Power balance at all buses is maintained by \eqref{eq:p_balance}--\eqref{eq:q_balance}.
Note that the formulation allows for multiple loads and generators with different individual limits at one bus.
In power flow equations \eqref{eq:line_p}--\eqref{eq:line_q}, power flows on each branch $(i,j)$ are defined. Note that each branch has two power flow equations for each active and reactive power, one in direction $(i,j)$ and another one in direction $(j,i)$ to reflect transmission losses. 
Branch apparent-power limits are constrained in \eqref{eq:linelim_from}.
As in most literature regarding BNI problems, transformers are modeled as branches, assuming no tap changing and no phase shifting.
In the formulations for load and generation limits in \eqref{eq:demlim}--\eqref{eq:qlim_gen}, it is assumed that each load and each generator can be continuously controlled within their respective operational ranges, which is a common assumption in BNI modeling \cite{Zhao.2013, Abedi.2021, Bienstock.2010}. 
Since there are no binary lower-level variables to switch off generators or demands, the lower bound for demand and generation in \eqref{eq:demlim}--\eqref{eq:plim_gen} is assumed to be $0$ to avoid infeasibilities \cite{Bienstock.2010}.
Since $B_{ij} \gg B_{ij}^\text{S}$ holds, $B_{ij}^\text{S}=0$ is assumed for computational performance in Sec. \ref{sec:case_studies}, with negligible effects on the results.
\par

\subsection{Linearization of nonconvex terms and MILP reformulation} \label{sec:linearized}

The model presented in Sec. \ref{sec:nonconvex} is an NP-hard nonconvex bilevel mixed-integer nonlinear program (BMINLP) due to its nonconvex constraints in \eqref{eq:line_p}--\eqref{eq:line_q}.
Several complexity reduction strategies exist, including linearization, linear relaxation, and convexification, as outlined in Sec. \ref{sec:introduction}. 
The most common relaxation is the DC approximation, applying strong assumptions about impedances, voltage magnitudes, and voltage angles, resulting in a fully linear lower level.
In this work, we adapt the linearizations from \cite{Abedi.2021} as described below, and use a DC model for comparison.

The nonconvex power flow equations \eqref{eq:line_p}--\eqref{eq:line_q} are linearized using a second-order Taylor linearization and a piecewise-linear approximation of quadratic terms, without introducing additional binary variables. 
This omission is only valid if branch losses are explicitly or implicitly penalized in the objective function.
For the presented models, increasing branch losses cannot improve the lower-level objective value since the minimum generation $\underline P_g\genmin=0$ at each generator. 
There might be cases in which branch losses have neither a positive nor a negative impact on the objective value.
In these cases, multiple feasible solutions may exist, potentially affecting generation values but not the attack vector or objective value.
Branch thermal limits are linearized using an $n$-sided inner polygon approximation ($n=8$), ensuring feasibility \cite{Akbari.2016}.
The resulting model is a bilevel mixed-integer bilinear program (BMIBLP), with bilinear terms arising from the products of binary upper-level variables $z_{ij}$ and continuous lower-level variables.

The DC OPF model is given by replacing the lower-level model in \eqref{eq:demlim}--\eqref{eq:linelim_from} with \eqref{eq:dc_demlim}--\eqref{eq:dc_lim_line} and is already a BMIBLP with decision variables $\Omega_{\text{DC}} \coloneqq \big\{p\dem_d, p\gen_g, \theta_i, p_{ij}^\text{K}, \big\} $.

\begin{subequations}\label{eq:dc_bilevel_problem}
\renewcommand{\theequation}{\theparentequation.\arabic{equation}} \vspace{-0.5em}
\begin{flalign}
  & 0 \leq p_d\dem \leq P_d\dem \quad \forall d \in D \label{eq:dc_demlim}\\
  & 0 \leq p_g\dem \leq P_g\gen \quad \forall g \in G \label{eq:dc_genlim}\\
  & p_{ij} = -z_{ij} B_{ij} (\theta_i - \theta_j)\quad \forall (i,j) \in K \label{eq:dc_line_flow}\\
  & -\overline S_{ij} \leq p_{ij} \leq \overline S_{ij} \quad \forall (i,j) \in K \label{eq:dc_lim_line}
\end{flalign}
\end{subequations}

All remaining bilinear terms are linearized using Big-M linearization following \cite{FortunyAmat.1981, Abedi.2021}, yielding a bilevel mixed-integer linear program (BMILP) in each case.

The BMILPs for the linearized AC (LAC) and DC OPF models are further reformulated as single-level mixed-integer linear programs (MILPs), allowing the use of general-purpose solvers. 
This reformulation applies strong duality to the lower-level problem, replacing it with its primal and dual constraints, along with the strong duality condition to define optimality.
Bilinear terms involving upper-level and dual lower-level variables are again linearized applying Big-M linearization, and the resulting problem is a MILP.

\section{New evaluation and scoring methodologies} \label{sec:methodologies}
This section presents two methodologies to enhance bilevel vulnerability assessments.
The first introduces an evaluation procedure for comparing different power flow formulations in the lower-level OPF, applied here to the LAC and DC models. 
The second applies CAV scoring across multiple load and generation cases.

%
Let $x$ denote the tuple of attacked branches out of the set $\mathcal{X}$ of all possible attack combinations, i.e., $x \in \mathcal{X}$ is one possible CAV. For each approach $a \in \left\{ \text{DC}, \text{LAC}\right\}$ and attack budget $\overline Z$, \xatt~is the $n^\text{th}$-worst attack combination, as defined in~\eqref{eq:xatt_def}, where $n$ defines the rank of the solution. 
This tuple contains the indices of branches $(i,j)$ for which $z_{ij}=0$ in an optimal solution (e.g., the tuple of attacked branches $x_{\mathrm{a},2}^{(1)} = \{(1,4), (2,5)\}$). 
Index $m(n)$ is the rank $n^\text{th}$-worst LAC attack vector in the set of DC solutions, as defined in \eqref{eq:rank_def}. 
\par \vspace{-0.5em}
\begin{subequations}\label{eq:analysis_defs}
\renewcommand{\theequation}{\theparentequation.\arabic{equation}} \vspace{-0.5em}
\begin{flalign}
  & \xattm \coloneqq \Big \{ (i,j) \in K \mid z_{ij}^{\star,a,\overline Z}=0 \Big \} \quad \forall n \in \mathcal{N}^{a} \label{eq:xatt_def}\\
  & m(n) = \operatorname{rank}_{\mathrm{DC},\overline Z}\!\left(x^{(n)}_{\mathrm{AC},\overline Z}\right) \quad \forall n \in \mathcal{N}^\text{LAC} \label{eq:rank_def}
\end{flalign}
\end{subequations}
Identifying only the optimal solution of a model is insufficient for systematic HILP event analysis, as multiple other CAVs may still exhibit high potential damage.
To address this, we adapt the procedure from \cite{Toenges.2025}, adding a constraint for each discovered CAV to exclude it from subsequent searches. 
Unlike the similar approach described in \cite{Ding.2017}, unattacked components are not explicitly considered.
If an optimal solution \xatt contains fewer than $\overline Z$ components, any possible \xattnew \ containing \xatt \ is excluded.
This ensures that only distinct, critical vulnerabilities are identified, as adding further components to \xatt \, does not improve the objective value $\obja$.
The total number of identified CAVs with approach $a$ is denoted as $N^a$, while $\mathcal{L}^\mathrm{a}$ is the set of all identified CAVs.

\subsection{Evaluation procedure for different OPF formulations} \label{sec:comp_method}
Previous studies have shown that optimal solutions and attack vectors can differ between AC and DC formulations \cite{Abedi.2021, LopezLezama.2017}.
What remains to be clarified is whether these differences stem from the DC approach completely failing to identify certain CAVs (e.g. when reactive-power constraints are binding in the AC model), or whether the DC model instead yields a different, yet still comparably critical, ordering of attack vectors within acceptable error margins.
To determine whether the DC approach fails to detect impactful CAVs or simply produces a different ranking, we propose the evaluation procedure outlined in Alg.~\ref{alg:acdc_comp}.
The procedure takes as input CAV lists generated by the LAC and DC models, denoted $\mathcal{L}^{\text{LAC}}$ and $\mathcal{L}^{\text{DC}}$. 
This procedure is not limited to the comparison here, but can also be applied without modification to evaluate other vulnerability assessment approaches.
\par
\begin{algorithm}[tb]
\footnotesize
\captionsetup{font=footnotesize} 
\caption{Evaluation of differences between OPF formulations}\label{alg:acdc_comp}
\begin{algorithmic}
\State get $\mathcal{L}^{\text{LAC}}, \mathcal{L}^{\text{DC}}$, initialize $n=1, U=0$ \Comment{Initialize analysis}
\While{$n \leq N^{\text{LAC}}$} \Comment{All considered LAC solutions} \vspace{0.2em}
\If{$\xattac \in \mathcal{L}^{\text{DC}}$} \Comment{Check if LAC attack also in DC results}
\State $\Delta \objabs = \objacm - \objdcm $ \Comment{Calculate absolute ...}
\State $\Delta \objrel = \tfrac{\Delta \objabs}{\objacm} $ \Comment{... and relative objective value gap}
\Else
\State $U \gets U+1$ \Comment{Increase number of undetected solutions}
\EndIf
\State $n \gets n+1$ \Comment{Increase solution index}
\State Store $\Delta \objabs$ and $\Delta \objrel$ in $\mathcal{L}^{\text{LAC}}_n$ \Comment{Save results}
\EndWhile \\
\Return $\mathcal{L}^{\text{LAC}}$
\end{algorithmic}
\end{algorithm}
\captionsetup{font=footnotesize}
Afterwards, we compute key performance indicators (KPIs) across all solutions, including the percentage of undetected CAVs \eqref{eq:acdc_hidden}, and the absolute and relative average deviations of $\obja$~\eqref{eq:acdc_offset_abs}--\eqref{eq:acdc_offset_rel}.
\begin{subequations}\label{eq:acdc_kpis}
\renewcommand{\theequation}{\theparentequation.\arabic{equation}}
\begin{flalign}
  & u = \tfrac{\hidden}{|\acset|} \label{eq:acdc_hidden} \\
  & \Psi^\text{abs} = \textstyle \sum_{\sol \in \acset}(\Delta \objabs) \cdot \tfrac{1}{|\acset|} \label{eq:acdc_offset_abs} \\
  & \Psi^\text{rel} = \textstyle \sum_{\sol \in \acset}(\Delta \objrel) \cdot \tfrac{1}{|\acset|} \label{eq:acdc_offset_rel}
\end{flalign}
\end{subequations}
\subsection{Scoring methodology for multiple load\,/\,generation cases} \label{sec:scoring_method}
Optimal attack vectors can vary depending on the specific load and generation configurations at different time steps.
While the worst-case solution for each time step is identified in \cite{Abedi.2020}, attack vectors that are near-optimal in one case may remain critical in others, whereas some CAVs depend more strongly on the specific load case. 
For grid operators, it is therefore important to identify which CAVs consistently pose significant risks across a range of scenarios, such as different time steps or grid configurations.
To this end, we introduce a multi-load-case scoring methodology in Alg.~\ref{alg:ts_scoring}.
This approach evaluates each CAV across all considered cases and produces two rankings: a rank score $\Phi^\text{rank}$ (sorted in ascending order), and an objective score $\Phi^\text{obj}$ (sorted in descending order).
\setlength{\textfloatsep}{1pt}
\begin{algorithm}[tb]
\footnotesize
\captionsetup{font=footnotesize} 
\caption{CAV scoring across multiple time steps}\label{alg:ts_scoring}
\begin{algorithmic}
\State get $T$ and all relevant simulation results for given $\overline Z$, initialize $t=1$
\State create $\mathcal{C}, \mathcal{R}, \mathcal{Y}$ \Comment{Each containing all possible attack combinations}
\While{$t \leq T$} \Comment{Consider all time steps}
\State $n=1$, get $\mathcal{L}^a(t)$ \Comment{Set solution index to 1, get solutions}
\While{$n \leq |\mathcal{L^a}(t)|$} \Comment{Consider all available solutions} \vspace{0.2em}
\State $\mathcal{C}\left(\xattm(t)\right ) \gets \mathcal{C}\left(\xattm(t)\right )+1$ \Comment{Update appearance counter} \vspace{0.2em}
\State $\mathcal{R}\left(\xattm(t)\right ) \gets \mathcal{R}\left(\xattm(t)\right )+n$ \Comment{Update rank sum} \vspace{0.2em}
\State $\mathcal{Y}\left(\xattm(t)\right ) \gets \mathcal{Y}\left(\xattm(t)\right )+\obja(t)$ \Comment{Update objective sum} \vspace{0.2em}
\State \Comment{$\mathcal{C}, \mathcal{R}, \mathcal{Y}$ are updated for the considered optimal attack vector}
\State $n \gets n+1$ \Comment{Increase solution index}
\EndWhile
\State $t \gets t+1$ \Comment{Go to next time step}
\EndWhile
\State $\Phi^\text{rank}(x) = \tfrac{\mathcal{R}(x)\cdot T}{(\mathcal{C}(x))^2} \; \forall x$ \Comment{Calculate rank score and ...} \vspace{0.2em} 
\State $\Phi^\text{obj}(x) = \tfrac{\mathcal{Y}(x)}{\mathcal{C}(x)}\cdot\tfrac{\mathcal{C}(x)}{T} \; \forall e$ \Comment{... objective score for each possible attack} \vspace{0.2em}
\Return $\Phi^\text{rank}, \Phi^\text{obj}$ \Comment{Return scoring results}
\end{algorithmic}
\end{algorithm}
\captionsetup{font=footnotesize} 

The objective score of a CAV, $\Phi^\text{obj}(x)$, is defined as its mean objective value across all time steps. The rank score, $\Phi^\text{rank}(x)$, corresponds to the average rank of a CAV in the time steps where it appears, scaled by the inverse of its occurrence frequency. 
For instance, a CAV that appears every time step with rank 3 receives the same rank score as one that appears every third time step but is consistently ranked first.
Because the scoring algorithm requires only the CAVs and their associated lost load, it can be applied broadly to any vulnerability assessment approach that provides this information.


\section{Case Studies} \label{sec:case_studies}
\subsection{Case study description}\label{sec:case_study_description}
All case studies are based on the open-source SimBench high-voltage test grid \verb|1-HV-urban--0-no_sw|~\cite{Meinecke.2020}, consisting of 82 buses, 79 loads, 98 generation units, 113 lines, 14 substations, three transformers, and one external grid connection. 
The dataset includes representative one-year time series for different types of load and generation at a 15-minute resolution for the reference year 2016 \cite{Meinecke.2020}. 
Throughout this section, the terms \textit{time step} and \textit{load and generation case} are used interchangeably to denote the discrete operating points considered.
\par
For the analyses in this section, two reference days are selected: January 29, with the largest positive difference between total apparent power load and generation, and May 30, with the largest negative difference. 
For the evaluation of the DC versus LAC approaches in Sec.~\ref{sec:comp_results}, eight time steps are taken from each day at three-hour intervals, yielding 16 time steps in total.
For the scoring methodology demonstration in Sec.~\ref{sec:scoring_results}, all time steps from each reference day together with the three preceding and following days are analyzed, yielding 1344 time steps in total. 
Note that the methodology itself is independent of the chosen model formulation.
Here, the primal-dual DC model from Sec.~\ref{sec:formulations} is used for illustration.

\subsection{Results for evaluating OPF formulations} \label{sec:comp_results}
\begin{figure}[b]
\centering
\includegraphics[width=8.2cm]{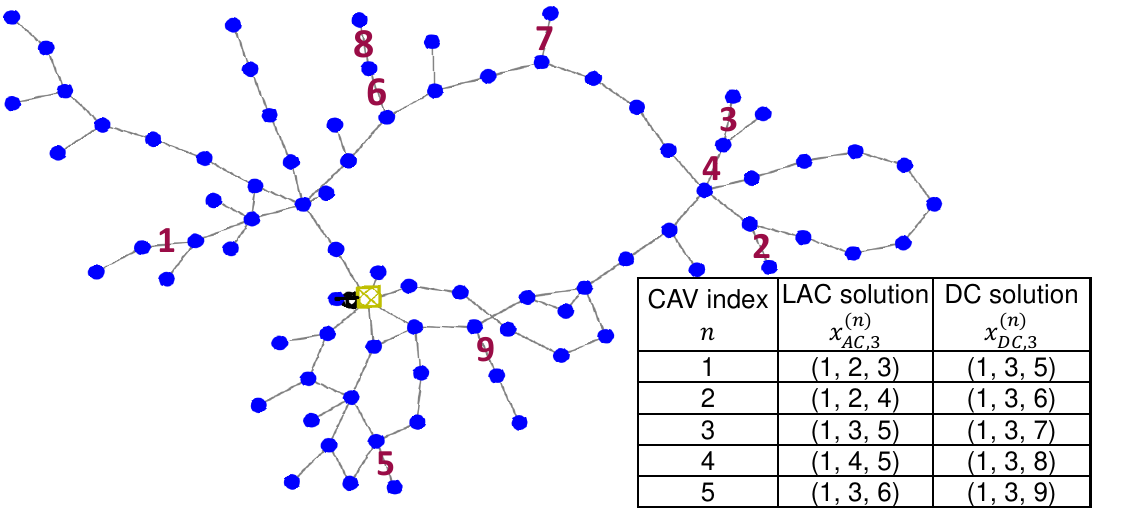}
\caption{SimBench high-voltage test grid topology with the best five CAVs for the DC and LAC approach (red numbers indicate attacked line indices).}
\label{fig:hv8_grid}
\end{figure}

The methodology from Sec.~\ref{sec:comp_method} is applied to compare the DC and LAC approaches introduced in Sec.~\ref{sec:formulations}. 
Fig.~\ref{fig:hv8_grid} shows the SimBench grid and highlights the top five CAVs identified by both models ($N^\text{LAC} = N^{\text{DC}} = 5$) at the time step May 30, 00:00, with $\overline Z=3$.
The results indicate that attacks predominantly target radial lines, which is expected given the high number of distributed generators and the fact that disconnecting a radially connected bus requires fewer attacks than isolating a bus in a meshed area. 
Although this observation is specific to the analyzed grid structure, the methodology itself is applicable to any topology, including fully meshed grids.
Notably, the optimal attack vectors differ: the DC model captures only two of the LAC solutions. 
This suggests that the DC approach may overlook relevant CAVs, a hypothesis that is further examined using Alg.~\ref{alg:acdc_comp}.
\par
\begin{figure}[t]
\centering
\includegraphics[width=8.83cm]{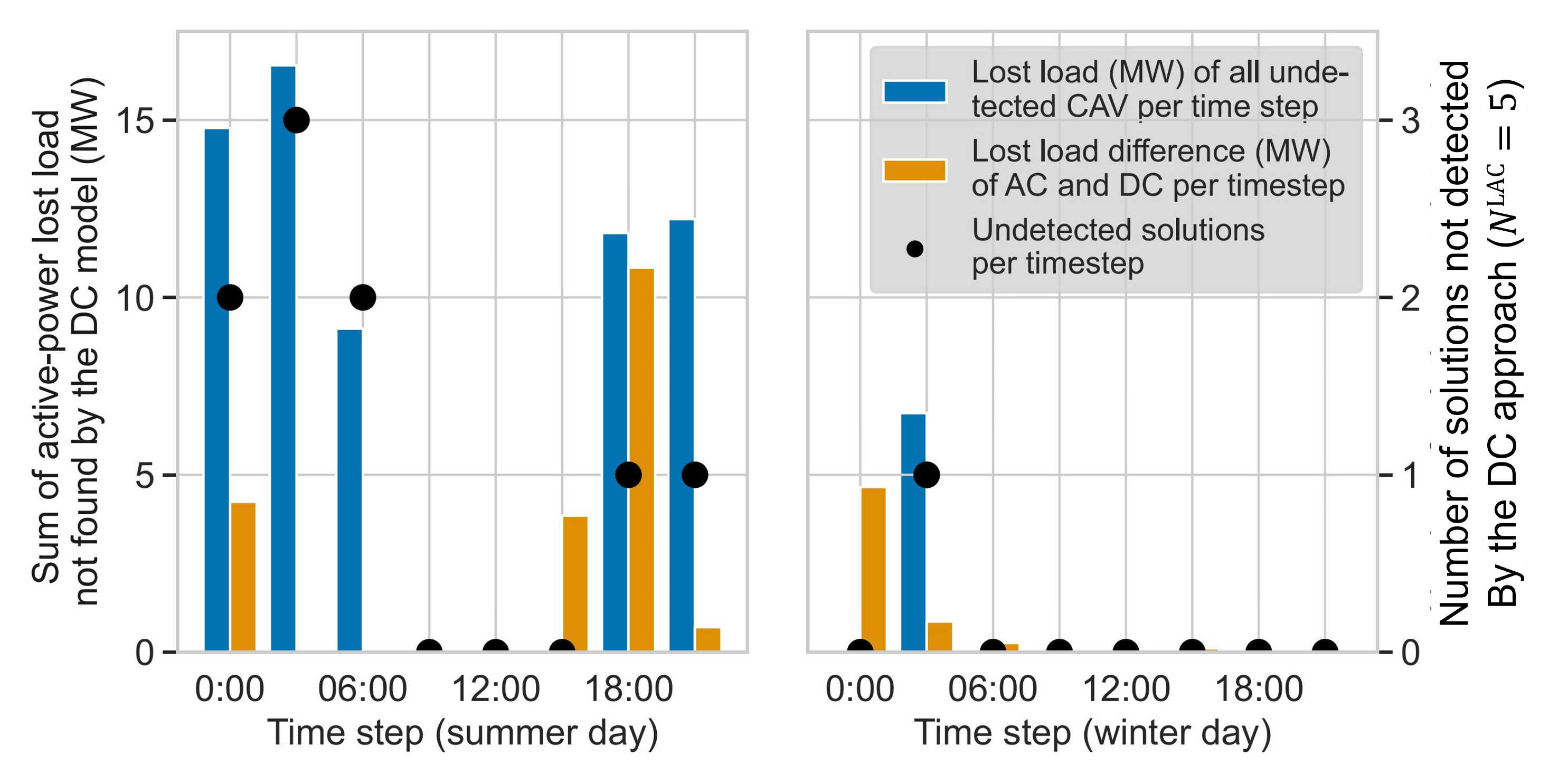} 
\caption{Sum of lost load $\zeta^{(n)}_{\text{LAC},2}$ for LAC CAVs that remain undetected with the DC approach (blue bars), sum of lost load underestimated with the DC approach for identified CAVs $\Delta \zeta^{(n)}_{\text{abs},2}$ (orange bars) and absolute number of undetected LAC CAVs per time step (dots).}
\label{fig:acdc_plot_16ts}
\end{figure}

Fig.~\ref{fig:acdc_plot_16ts} investigates this issue for $\overline Z=2$ across all 16 selected time steps. For each time step, five LAC solutions are computed, while at least 50 DC solutions are considered (or more if the 50\textsuperscript{th} DC solution still exceeds 50\,\% of the worst-case lost load). 
For instance, in the first summer time step, two of the five worst-case LAC CAVs are not detected by the DC model, corresponding to nearly 15\,MW of lost load. 
For three jointly identified CAVs, the DC model underestimates lost load by about 4\,MW. 
\par 
Undetected solutions occur on both days but are more frequent during the summer day, especially in the morning and evening, coinciding with periods of lower renewable generation.
The main insights are as follows.
First, the DC approach exhibits significant blind spots that extend beyond tolerable inaccuracies.
In the present case study, it fails to capture more than 15\,\% of the LAC solutions, with the corresponding KPIs $u=15.625\,\%$, $\Psi^\text{abs}=0.4748$\,MW, and $\Psi^\text{rel}=0.0432$. 
This implies that grid operators relying solely on the DC model may remain unaware of critical vulnerabilities.
Comparable results are observed for $\overline Z=1$ and $\overline Z=3$. 
These findings confirm and extend the previous results ~\cite{LopezLezama.2017, Abedi.2021}, showing that the discrepancies go beyond worst-case inaccuracies.
Second, the substantial variation across time steps highlights the strong influence of load and generation patterns, motivating the analysis in Sec.~\ref{sec:scoring_results}.
Although demonstrated here with LAC and DC models, the evaluation methodology is not restricted to these formulations. 
It extends prior worst-case-focused evaluations by quantifying systematic differences in vulnerability assessments across various approaches and can be applied to any grid topology or size, including real-world systems.

\subsection{Results for scoring multiple load\,/\,generation cases}\label{sec:scoring_results}
For all 1344 time steps described in Sec.~\ref{sec:case_study_description}, CAVs with lost load of at least 50\,\% of the worst case per time step are computed for each $\overline Z \in \{1, 2, 3\}$, resulting in $6{,}448$, $51{,}874$, and $41{,}217$ solutions, respectively. 
The objective scores $\Phi^\text{obj}$ and rank scores $\Phi^\text{rank}$ are then determined with Alg.~\ref{alg:ts_scoring}.
\par
\begin{figure}[tb]
\centering
\includegraphics[width=8.83cm]{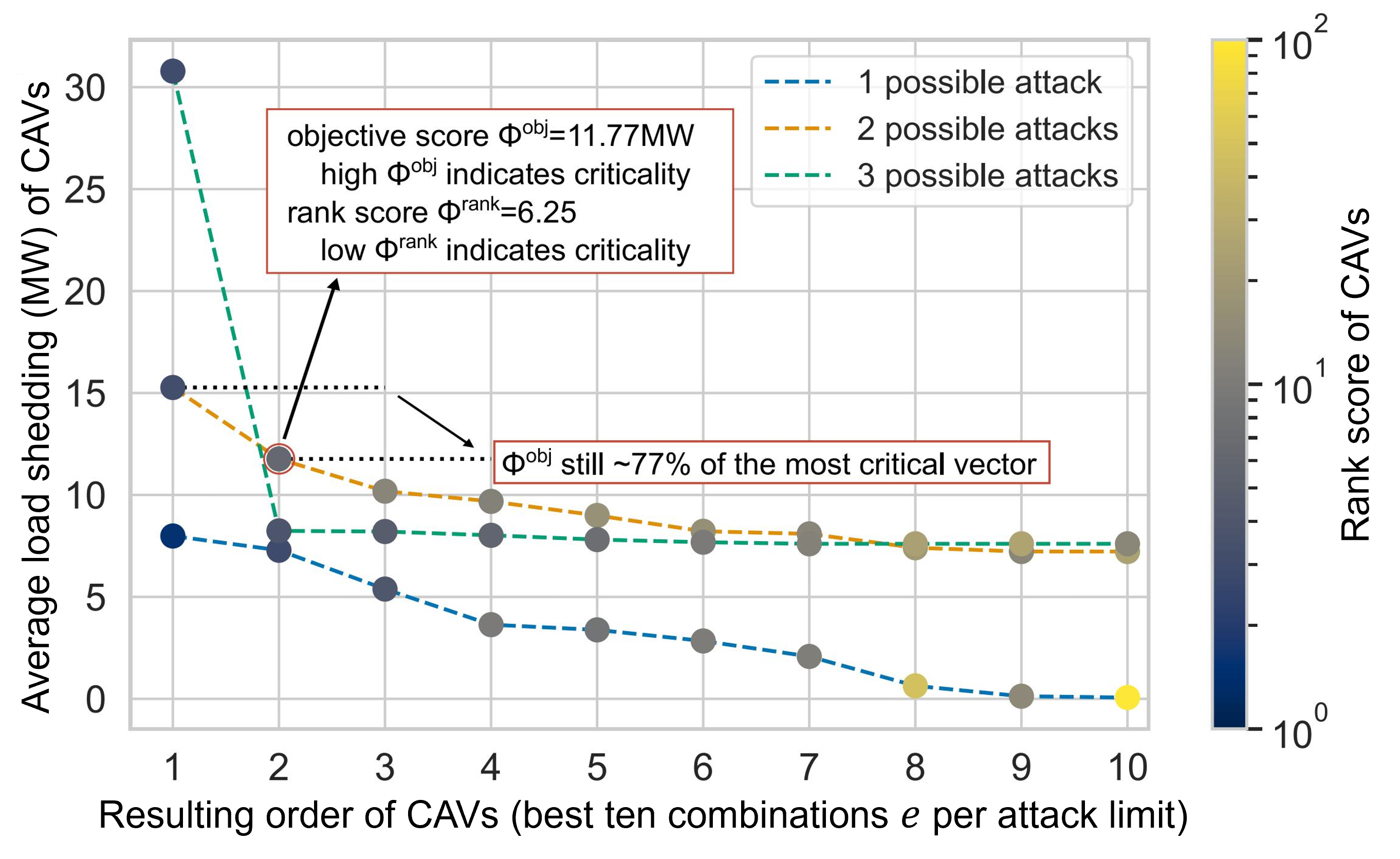} 
\caption{Objective and rank score of best 10 attack vectors for 1, 2 and 3 possible attacks, indexed by the $1^{\text{st}}$ to $10^{\text{th}}$ best combination per $\overline Z$.}
\label{fig:scoring_res}
\end{figure}

Fig.~\ref{fig:scoring_res} shows the ten CAVs with the highest objective scores for each $\overline Z$. 
For example, index 5 corresponds to the CAV with the fifth-highest objective score for each $\overline Z$. 
The corresponding rank score for each CAV is indicated by color. 
While objective and rank scores often align, several notable exceptions occur.
The objective score highlights CAVs that consistently cause high lost load and, therefore, require prioritization.
In contrast, the rank score captures how frequently CAVs appear and whether they are among the top solutions.
The combination of both measures provides particularly valuable insights.
For example, the CAV with a rank score of $6.25$ but solution index $2$ in Fig.~\ref{fig:scoring_res} would often remain undetected without cross-time-scoring, despite accounting for more than 77\,\% of the maximum average lost load. 
For $\overline Z=3$, a notable drop occurs between the highest and second-highest objective scores, which is explained by the top solution involving attacks on all three transformers.
By contrast, the remaining objective scores in Fig.~\ref{fig:scoring_res} decline more gradually with increasing index.
The obtained results underscore the importance of considering multiple CAVs beyond the worst case for a comprehensive vulnerability assessment. 
The proposed scoring methodology can be combined with existing approaches \cite{Abedi.2020, Sundar.2024} to provide a broader view of system vulnerabilities and to support resilience enhancement.

\section{Conclusion} \label{sec:conclusion}
This work introduces two methodologies to enhance vulnerability assessment in power systems: an evaluation procedure for BNI models with different OPF formulations, and a scoring methodology to analyze critical attack vectors (CAVs) across multiple load and generation scenarios. 
Both approaches move beyond traditional worst-case analysis to provide a more comprehensive understanding of system vulnerabilities.
\par
The comparative analysis of linearized AC (LAC) and DC OPF formulations demonstrates that the DC approach fails to identify a substantial share of critical vulnerabilities detected by the LAC model, even when multiple solutions are considered. 
This exposes the risk of relying solely on DC approximations in practice. 
At the same time, the computational burden of the AC formulation, even in its linearized form, highlights the ongoing need for efficient and scalable approximations in BNI modeling.
\par
Applying the scoring methodology across time steps reveals the strong dependency of CAVs on load and generation conditions.
This demonstrates the necessity of considering a range of grid configurations beyond worst cases to uncover vulnerabilities that might otherwise remain hidden.
\par
Overall, the proposed methodologies enhance the identification and understanding of HILP events, a key step toward developing robust detection and defense measures against adverse events.
While optimal defense strategies have been explored in earlier work, future research should focus on integrating them with advanced vulnerability assessment methods and with emerging cyber-physical energy systems to strengthen power system resilience.

\bibliographystyle{IEEEtran}
\bibliography{literature}

\end{document}